\documentclass[a4paper,12pt]{article}

\usepackage{amsfonts}
\usepackage{amsthm}
\usepackage{amsmath}
\usepackage{amsfonts}
\usepackage{latexsym}
\usepackage{amssymb}

\newtheorem{fed}{Definition}[section]
\newtheorem{teo}[fed]{Theorem}
\newtheorem*{teo*}{Theorem}
\newtheorem{lem}[fed]{Lemma}
\newtheorem{cor}[fed]{Corollary}
\newtheorem{pro}[fed]{Proposition}
\theoremstyle{definition}
\newtheorem{rem}[fed]{Remark}
\newtheorem{rems}[fed]{Remarks}
\newtheorem{exa}[fed]{Example}
\newtheorem{exas}[fed]{Examples}
\newtheorem{num}[fed]{}

\def\noi{\noindent}
\def\QED{\hfill $\blacksquare$}
\def\EOE{\hfill $\blacktriangle$}

\def\fii{\varphi }
\def\la{\lambda}

\def\cA{\mathcal{A}}
\def\cB{\mathcal{B}}
\def\cC{\mathcal{C}}

\def\cE{\mathcal{E}}
\def\cH{\mathcal{H}}
\def\cK{\mathcal{K}}

\def\cZ{\mathcal{Z}}

\def\casa{\mathcal{A}_{sa}}

\def\cstarr{$C^*$-algebra}
\def\cstar{$C^*$-algebra $\;$}
\def\cstars{$C^*$-algebras $\;$}

 \DeclareMathOperator{\tr}{tr}

\DeclareMathOperator*{\inversible}{GL}

    \DeclareMathOperator{\leqrt}{\precsim}

    \DeclareMathOperator{\leqm}{\preccurlyeq}

\newcommand{\pint}[1]{\displaystyle \left \langle #1 \right\rangle}

\newcommand{\hil}{\mathcal{H}}
\newcommand{\op}{L(\mathcal{H})}
\newcommand{\opsa}{L_{sa}(\mathcal{H})}
\newcommand{\posop}{L(\mathcal{H})^+}

\newcommand{\mat}{\mathcal{M}_n}
\newcommand{\matsa}{\mathcal{M}_n^{sa}}
\newcommand{\matpos}{\mathcal{M}_n^+}

\newcommand{\spec}[1]{\sigma\left( #1\right)}
\newcommand{\avi}[2]{\la_{#1}\left( #2\right)}

\newcommand{\ninf}[1]{\left\| #1 \right\|_{\infty}}

\newcommand{\espro}[2]{E_{#2}\left[#1\right]}

\begin{document}


\title{{\Huge\textbf{Jensen's inequality and majorization.}}}
\author{Jorge Antezana, Pedro Massey and Demetrio Stojanoff}
\date{}
\maketitle
\centerline{\it Dedicated to the memory of Gert K. Pedersen}

\bigskip

\begin{abstract}
Let $\mathcal{A}$ be a $C^*$-algebra and $\phi:\cA\rightarrow L(H)$ be a
positive unital map. Then, for a convex function $f:I\rightarrow
\mathbb{R}$ defined on some open interval and a self-adjoint
element $a\in \mathcal{A}$ whose spectrum lies in $I$, we obtain a
Jensen's-type inequality $f(\phi(a)) \ \leq \ \phi(f(a))$ where $\le$
denotes an operator preorder (usual order, spectral preorder,
majorization) and depends on the class of convex functions
considered i.e., operator convex, monotone convex and arbitrary
convex functions. Some extensions of Jensen's-type inequalities to
the multi-variable case are considered. 
\end{abstract}

\vglue1truecm

\noi {\bf Keywords:} Jensen's inequality, convex functions,
positive maps, majorization, spectral preorder.

\medskip
\noi {\bf AMS Subject Classifications:} Primary 47A63 ; Secondary 46L05, 46L10, 47A60, 15A45.



\section{Introduction}

Jensen's inequality is the continuous version of the usual
definition of convex function and it can be stated in the
following way: let $(X,P)$ be a probability space and let $g$ be a
bounded function such that for all $x\in X$, $g(x)\in I$. Then
\[
f\left(\int_Xg\;dP\right)\leq \int_Xf\circ g\;dP.
\]

In the context of \cstars the simplest generalization of Jensen's
inequality can be made by taking a state $\phi$ and a selfadjoint
element $a$ of a \cstar $\cA$ such that the spectrum of $a$ is
contained in $I$, then it holds that
\begin{equation}\label{jensen con estado en la introduccion}
f(\fii(a))\leq \fii(f(a)).
\end{equation}

The same inequality, replacing the state $\fii$ by a positive
unital map $\phi:\cA\rightarrow \cB$ between two \cstars
($\phi(a)\geq 0$ whenever $a\geq 0$ and $\phi(I)=I$) is false
unless $f$ be an operator convex function. Indeed T. Ando proved
that, given a Hilbert space $\hil$ and a function $F$ defined on some
open interval $I$, the following statements are equivalent:
\begin{description}
  \item[1.] $F$ is operator convex.
  \item[2.] $F(v^*a v)\leq v^* F(a)v$ for every $a\in L_{sa}(\cK)$ such that
  $\sigma(a)\subseteq I$, where $\cK$ is some Hilbert space and $v$ is any isometry
  from $\hil$ into $\cK$.
  \item[3.] $F(\phi(a))\leq \phi(F(a))$  for every $a\in\opsa$ such that
  $\sigma(a)\subseteq I$, and for every unital completely positive map $\phi:\op\rightarrow\op$.
  \item[4.] $F(\phi(a))\leq \phi(F(a))$ for every $a\in\opsa$ such that
  $\sigma(a)\subseteq I$, and for every unital positive map $\phi:\op\rightarrow\op$.
\end{description}

Although in some particular cases (e.g. when $\phi(a)$ and
$\phi(f(a))$ commute) we can obtain Jensen's type inequalities
like in item 4 for arbitrary convex functions, Ando's result tells
us that in general we can not consider the usual order for these
kind of inequalities and convex funtions which are not operator
convex.

However, previous works on the matter, such as Brown-Kosaki's and
Hansen-Pedersen's (\cite{[BrKo]}, \cite{[HanPed]}), suggest the
idea of considering Jensen's type inequalities with respect to
other preorders, such as the spectral and majorization preorders
(see subsection \ref{preordenes en algebras de operadores} for
their definitions).

In this paper we study different Jensen's-type inequalities in
which, as in the case of Ando's Theorem, the roll of
non-commutative integral is played by a positive (unital) map.

The paper is organized as follows. In section 2 we recall some
definitions and we introduce notation used throughout this paper.
For the sake of completness, we sketch the proof of Ando's Theorem
and give a short proof of Jensen's inequality for states in a
$C^*$-algebra.

 Section 3 is divided in two subsections, depending
on the preorder relation and the assumptions on the convex
functions.
 In subsection 3.1 we consider monotone convex functions and,
following Brown-Kosaki's ideas on the matters, we prove Jensen's
type inequalities with respect to the spectral preorder.
In subsection 3.2 we obtain Jensen's type inequalities for
arbitrary convex functions either by taking restrictions on the
algebra $\cB$ or by using the majorization preorder.

We collect the results of this section in the following list:
Let $\cA$, $\cB$  be unital $C^*$-algebras,
$\phi:\cA\rightarrow\cB$ a positive unital map, $f$ a
convex function defined on an open interval $I$ and $a\in \cA$, such that
$a = a^*$ and $\sigma(a) \subseteq I$. Then
\begin{enumerate}
\item If $f$ is monotone and $\cB$ is a von Neumann algebra, then
$$ f(\phi(a))\leqrt \phi(f(a))\,\,\, \text{(spectral preorder)}.$$

\item If $\cB$ is abelian or, more generally, if $\phi(f(a))$ and $\phi(a)$ commute,
then
$$f(\phi(a))\le \phi(f(a)).$$
\item If $E : \cB \to \cC$ is a conditional expectation and $\phi(a)$ belongs to the
centralizer of $E$, then $$\cE(g[f(\phi(a))])\leq
\cE(g[\phi(f(a))])$$  for every continuous increasing convex
function $g$ such that Im$(f)\subseteq$Dom$(g)$.

\item if  $\cB$ a finite factor, then
$$f(\phi(a)))\leqm \phi(f(a))\,\,\,\text{(majorization)}.$$
\end{enumerate}
We remark that some of this inequalities still hold for
$contractive$ positive maps, under the assumption that $0\in I$
and $f(0)\leq 0$.

We should also mention that in the item 3 of the previous list, a
similar inequality without the convex function $g$ can be obtain
as a corollary of Hansen and Pedersen's results \cite{[HanPed]}
and Stinespring's
representation theorem. However, this convex function is
meaningful because it let us connect Hansen and Pedersen's work
with the theory of majorization in finite factors.

 In section 4 we describe briefly the multi-variable
functional calculus and obtain in this context similar results to
those of section 3.2 by using essentially the same techniques.

 In section 5 we apply the results obtained to the finite
dimensional case, where there exist fairly simple expressions for
the spectral preorder  and for the majorization; we show that our
results generalize those appeared in a recent work by J. S. Aujla
and F. C. Silva \cite{[AujSil]}.

 Finally, section 6 is devoted to obtain some inequalities by
choosing particular functions and (unital) positive maps. For
instance, we prove non-commutative versions of the information
inequality, Liapunov's inequality and H\"older's inequality.

\section{Preliminaries}

Let $\cA$ be a $C^*$-algebra, throughout this paper $\casa$ will
be the real vector space of self-adjoint elements of $\cA$,
$\cA^+$ the cone of positive elements, $\inversible(A)$ the group
of invertible elements and $M_n(\cA)$ the $n\times n$ matrices
whose entries are elements of $\cA$. We also suppose that all the
\cstars in consideration are unital. In a similar way given a
Hilbert space $\hil$, $\op$ will be the algebra of all linear
bounded operators on $\hil$, $\opsa$ the real vector space of
self-adjoint operators, and $\posop$ the cone
 of all positive operators. For every $C \in \op$ its range
will be denoted by $R(C)$, its null space will be denoted by
$\ker(C)$ and its spectrum by $\spec C$.

On the other hand, if $p$ and $q$ are orthogonal projections, the
orthogonal projection onto the intersection of their ranges will
be denoted $p\wedge q$ and the orthogonal projection onto the
closed subspace generated by their ranges will be denoted $p\vee
q$

\subsection{Positive and Completely Positive Maps.}

To begin with, let us recall some definitions.

\medskip

\begin{fed}\label{definicion de positivo}\rm
Let $\cA$ and $\cB$ two $C^*$-algebras and $\phi:\cA\rightarrow
\cB$ a linear map. We shall say that $\phi$ is \textbf{positive}
if it maps positive elements of $\cA$ to positive elements of
$\cB$.
\end{fed}

\smallskip

\begin{fed}\label{definicion de completamente positivo}\rm
Given a positive map, $\phi$, between              two \cstar $\cA$ and $\cB$,
let us call $\phi_n$ the map between $M_n(\cA)$ and $M_n(\cB)$
defined by $\phi_n((a_{ij}))=(\phi(a_{ij}))$. We say that
$\phi$ is \textbf{n-positive} if $\phi_n$ is positive and
\textbf{completely positive} if $\phi$ is n-positive for all
$n\in\mathbb{N}$.
\end{fed}

\bigskip

The following result due to Stinespring is one of the most general
theorems which characterizes completely positive maps from a
\cstar into $\op$.

\medskip

\begin{teo}\label{stinespring caraterizacion}
Let $\cA$ be a \cstar and let $\phi:\cA\rightarrow \op$ be a
completely positive map. Then, there exists a Hilbert space $\cK$,
a unital $^*$-homomorphism $\pi:\cA\rightarrow L(\cK)$, and a
bounded operator $V:\cH\rightarrow\cK$ with $\|\phi(1)\|=\|V\|^2$
such that
\[
\phi(a)=V^*\pi(a)V.
\]
\end{teo}

\medskip

The reader can find a proof of this theorem in \cite{[Pau]} or in
\cite{[Hosh]}. Note that in our case, since the map will be either
unital or a contraction, the bounded operator $V$ will be either
an isometry or a contraction.

Another important result about completely positive maps is the
next theorem, also due to Stinespring (See \cite{[Pau]}).

\smallskip

\begin{teo}\label{stinespring conmutativo}
Let $X$ be a compact Hausdorff space and $C(X)$ the \cstar of
continuous function on $X$. Then, every positive map defined on
$C(X)$ is completely positive.
\end{teo}

\subsection{Convex and Operator Convex Maps.}

\begin{fed}\label{convex}\rm
A real function $f$ defined on a segment $(a,b)$, where
$-\infty\leq a<b\leq\infty$, is called \textbf{convex} if the
inequality
\begin{equation}\label{ecuacion de convexa}
f(\la x+(1-\la) y)\leq \la f(x)+(1-\la)f(y)\, ,
\end{equation}
holds whenever $a<x<b$, $a<y<b$ and $0\leq \la\leq 1$.
A function $g$ such that $-g$ is convex is called concave.
\end{fed}

\medskip
\noindent
It is not difficult to prove that equation (\ref{ecuacion de convexa})
is equivalent to the requirement that
\begin{equation}\label{otra ecuacion de convexa}
\frac{f(t)-f(s)}{t-s}\leq \frac{f(u)-f(t)}{u-t}\, ,
\end{equation}
whenever $a<s<t<u<b$.

\begin{pro}\label{propiedad lla}
Let $f$ be a convex function defined in an open interval $I$ and
let $K\subseteq I$ be compact. Then, there exist a (countable)
family of linear functions $\{f_i\}_{i\geq 1}$ such that for each
$x\in K$ it holds that
\[
f(x)=\sup_{i\geq 1}f_i(x).
\]
\end{pro}

\medskip

\begin{fed}\label{operator convex}\rm
A continuous function $F:I\rightarrow \mathbb{R}$ is said to be
\textbf{operator convex} if
\begin{equation}\label{ecuacion de operator convex}
F(\la a+(1-\la) b)\leq \la F(a)+(1-\la)F(b) \, ,
\end{equation}
for each $\la\in [0,1]$ and every pair of self-adjoint operators
$a,b$ on an infinite dimensional Hilbert space $\hil$ with
spectrum in $I$. As in the scalar case, $F$ is called operator
concave if $-F$ is operator convex.
\end{fed}

\begin{rems}\label{obs conv de op}
$\;$
\begin{description}
  \item[1.] The equation (\ref{ecuacion de operator convex}) is
  well defined because the set of self-adjoint operator with
  spectrum in a segment $I$ is convex.
  \item[2.] The function $F$ is called \textbf{matrix convex of
  order n} if equation (\ref{ecuacion de operator convex}) is
  satisfied for operators belonging to an n dimensional Hilbert space.
  It is well known that a function is operator convex if and
  only if it is matrix convex of arbitrary orders. The interested
  reader can find a proof of this fact in \cite{[BenSher]}.
  \item[3.] Clearly, an operator convex function is convex but the
  converse is not true. Indeed it is known that the class of
operator convex functions has a very rich structure. For example
(see \cite{[Bat]}), a function $f:[0,\infty)\rightarrow
\mathbb{R}$ is an operator convex function if and only if there
exist a positive finite measure $\mu$ on $[0,\infty)$ and real
constants $\alpha,\, \beta,\,\gamma$ with $\gamma\geq0$, such that
$$
f(t)=\alpha+ \beta t+\gamma t^2 + \int _{0} ^{\infty}
\frac{\lambda t^2}{1+\lambda t}\, d\mu(\lambda)
$$
In particular $f$ is infinitely differentiable on $(0,\infty)$.
 \EOE
\end{description}
\end{rems}

Some well known functions that are operator convex and operator
concave are listed below:

\begin{exas}$\;$
\begin{description}
  \item[1.] $F(x)=x^r$ for $r\in[-1,0]\cup [1,2]$ is operator convex and
  for $r\in [0,1]$ is operator concave.
  \item[2.] $F(x)=\log(x)$ and $F(x)=- x\log(x)$ are operator concave in
  $(0,\infty)$. \EOE
\end{description}
\end{exas}

The notion of convexity is very close to Jensen's inequality. The
following Jensen-type inequality for states defined on a \cstar is
the simplest generalization to a non-commutative context.

\begin{teo}\label{jensen para estados}
Let $\cA$ be a \cstar, $\fii:\cA\rightarrow \mathbb{C}$ a state
and $f$ a convex function defined on some interval
$(\alpha,\beta)$ ($-\infty\leq \alpha<\beta\leq \infty$). Then it
holds that
\[
f(\fii(a))\leq \fii(f(a))) \, ,
\]
for every $a\in\casa$ whose spectrum is contained in
$(\alpha,\beta)$.
\end{teo}
\paragraph{Proof.}
Put $\tau=\fii(a)$. Then $\alpha<\tau<\beta$. If $\gamma$ is the
supremum of quotients on the left of (\ref{otra ecuacion de
convexa}), where $\alpha<s<\tau$, then $\gamma$ is no larger than
any of the quotients on the right of the same inequality, for
$\mu\in(\tau,\beta)$. Then
\[
f(\sigma)\geq
f(\tau)+\gamma(\sigma-\tau)\hspace{1cm}(\alpha<\sigma<\beta).
\]
Hence
\[
f(a)-f(\tau)+\gamma(a-\tau)\geq 0.
\]

Finally, if we apply $\fii$ to this inequality, taking into
account that $\fii$ is a state we get the desired inequality.\QED

\medskip
There also exist a closed relation between the notion of operator
convexity and Jensen-type inequalities, as the following theorem
due to Ando shows.

\begin{teo}\label{jensen pesadazo}
Let $\hil$ be an infinite dimensional Hilbert space and $F$ a function defined on some
open interval $I$.  Then, the following statements are
equivalent:
\begin{description}
  \item[1.] $F$ is operator convex.
  \item[2.] $F(v^*a v)\leq v^* F(a)v$ for every $a\in L_{sa}(\cK)$ such that
  $\sigma(a)\subseteq I$, where $\cK$ is some Hilbert space and $v$ is any isometry
  from $\hil$ into $\cK$.
  \item[3.] $F(\phi(a))\leq \phi(F(a))$  for every $a\in\opsa$ such that
  $\sigma(a)\subseteq I$, and for every unital completely positive map $\phi:\op\rightarrow\op$.
  \item[4.] $F(\phi(a))\leq \phi(F(a))$ for every $a\in\opsa$ such that
  $\sigma(a)\subseteq I$, and for every unital positive map $\phi:\op\rightarrow\op$.
\end{description}
\end{teo}
\paragraph{Sketch of proof.}
\begin{description}
  \item[$1. \Leftrightarrow 2.$] This is a well known fact. A proof
  can be found in \cite{[Bat]}.
  \item[$2.\Rightarrow 3.$] Given a unital completely positive map
  $\phi$, by theorem \ref{stinespring caraterizacion} there exists a Hilbert space
  $\cK$, a representation $\pi:\cA\rightarrow L(\cK)$, and a
  isometry $K:\cH\rightarrow\cK$ such that:
\[
\phi(a)=K^*\pi(a)K.
\]
So, as $F(\pi(a))=\pi(F(a))$, using (2.) we have that:
\[
F(\phi(a))=F(K^*\pi(a)K)\leq
K^*F(\pi(a))K=K^*\pi(F(a))K=\phi(F(a)).
\]
  \item[$3.\Rightarrow 4.$] The unital positive map restricted to
  the unital
  abelian \cstar generated by $a$ is completely positive according
  to the theorem \ref{stinespring conmutativo}. So (4.) is a
  consequence of (3.).
  \item[$4.\Rightarrow 2.$] It is clear. \QED
\end{description}

\subsection{Spectral Preorder and Majorization.}\label{preordenes en algebras de operadores}

In what follows $\espro{I}{a}=\chi_{I}(a)$ denotes the spectral
projection of a self-adjoint operator $a$ in a von Neumann algebra
$\cA$, corresponding to a (Borel) subset $I\subseteq \mathbb{R}$.

To begin with, let us recall the notion of spectral preorder.

\begin{fed}\label{definicion spectral preorder}\rm
Let $\cA$ be a von Neumann algebra.
Given $a,b\in \cA_{sa}$, we shall say that $a\leqrt b$ if
$\espro{(\alpha,+\infty)}{a}$ is equivalent, in the sense of
Murray-von Neumann, to a subprojection of
$\espro{(\alpha,+\infty)}{b}$ for every real number $\alpha$.
\end{fed}

In finite factors the following result can be proved.

\begin{pro}\label{proposicion para spectral preorder}
Let $\cA$ be a finite factor, with normalized
 trace $tr$. Given $a,b\in\cA_{sa}$, the
following conditions are equivalent:
\begin{enumerate}
  \item $a\leqrt b$.
  \item $\tr(f(a))\leq \tr(f(b))$ for every continuous increasing function $f$ defined
  on an interval containing both $\spec{a}$ and $\spec{b}$.
\end{enumerate}
\end{pro}

\noi
Finally, Eizaburo Kamei defined the notion of majorization and
submajorization in finite factors (see \cite{[Ka]}):

\begin{fed}\label{definicion mayorizacion}\rm
Let $\cA$ be a finite factor with normalized trace $tr$.
Given $a,b\in \cA_{sa}$, we shall say that $a$ is submajorized by
$b$, and denote $a\leqm b$, if the inequality:
\[
\int_{0} ^\alpha e_{a}(t)\,dt\leq\int_{0} ^\alpha e_{b}(t)\,dt
\]
holds for every real number $\alpha$, where
$e_{c}(t)=\inf\{\gamma:\, tr(E_{c}[(\gamma, \infty)])\leq t\}$ for
$c\in \cA_{sa}$. We shall say that $a$ is majorized by $b$ if
$a\leqm b$ and $\tr(a)=\tr(b)$.
\end{fed}

\noi
The following characterizarion of submajorization also appears in \cite{[Ka]}.
\begin{pro}\label{proposicion mayorizacion}
Let $\cA$ be a finite factor with normalized trace $tr$.
Given $a,b\in\cA_{sa}$, the following conditions are equivalent:
\begin{enumerate}
  \item $a\leqm b$.
  \item $\tr(f(a))\leq \tr(f(b))$ for every increasing convex function $f$ defined
  on an interval containing both $\spec{a}$ and $\spec{b}$.
\end{enumerate}
\end{pro}

\section{Jensen-type inequalities.}\label{mainsection}

Throughout this section $\phi$ is a positive unital map from a
$C^*$-algebra $\cA$ to another $C^*$-algebra $\cB$,
$f:I\rightarrow \mathbb{R}$ is a convex function defined on the
open interval $I$ and $a\in \cA_{sa}$ whose spectrum lies in $I$.
Note that the spectrum of $\phi(a)$ is also contained in $I$.

As we mentioned in the introduction, by Ando's Theorem, we can not
expect a Jensen-type inequality of the form
\begin{equation}\label{Jensen de Ando}
f(\phi(a))\leq \phi(f(a))
\end{equation}
for an arbitrary convex function $f$ without other assumptions.
This is the reason why, in order to study inequalities similar to
(\ref{Jensen de Ando}) for different subsets of convex function,
we shall use the spectral and majorization (pre)orders, or we
shall change the hypothesis made over the \cstar where the
positive map takes its values.

Although most of the inequalities considered in this section
involve unital positive maps, similar results can be obtained for
contractive positive maps by adding some extra hypothesis on  $f$.

\subsection{Monotone convex and concave functions. Spectral Preorder.}

\medskip

In this section we shall consider monotone convex and concave
functions. The following result due to Brown and Kosaki
(\cite{[BrKo]}) indicates that the appropriate order relation for
this class of functions is the spectral preorder. Let $\cA$ be a
semi-finite von Neumann algebra, and let $v\in\cA$ be a
contraction; then, for every positive operator $a\in\cA$ and every
continuous monotone convex function $f$ defined in $[0,+\infty)$
and such that $f(0)=0$, they proved that
\[
v^*f(a)v\leqrt f(v^*a v).
\]

The following Theorem is an analogue of Brown and Kosaki's result,
in terms of positive unital maps and monotone convex functions.
The proof we give below follows essentially the same lines.

\medskip
\begin{teo}\label{teo de jensen con el orden retruch}
Let $\cA$ be an unital $C^*$-algebra, $\cB$ a von Neumann algebra
and $\phi:\cA\rightarrow\cB$ a positive unital map. Then, for
every monotone convex function $f$, defined on some interval $I$,
and for every self-adjoint element $a\in\cA$ whose spectrum is
contained in $I$ it holds that
\begin{equation}\label{jensen orden retruch}
f(\phi(a))\leqrt \phi(f(a))
\end{equation}
\end{teo}
\paragraph{Proof.}

According to the definition, given $\alpha\in\mathbb{R}$ we have
to prove that there exists a projection $q\in \cA$ such that
\[
\espro{(\alpha,+\infty)}{f(\phi(a))}=\espro{\{f>\alpha\}}{\phi(a)}\sim
q \leq \espro{(\alpha,+\infty)}{\phi(f(a))}.
\]

To begin with, we claim that $\espro{\{f >
\alpha\}}{\phi(a)}\wedge \espro{(-\infty,\alpha]}{\phi(f(a))}=0$.
In fact, take $\overline{\eta}\in$ $ \displaystyle
R\left(\espro{\{f > \alpha\}}{\phi(a)}\right)$,
$\|\overline{\eta}\|=1$. Since $f$ is monotone we have that
\[
\alpha < f(\pint{\phi(a)\overline{\eta},\;\overline{\eta}})\, ,
\]
and using Theorem \ref{jensen para estados} we get
\[
\alpha < \pint{\phi(f(a))\overline{\eta},\;\overline{\eta}}.
\]
But, on the other hand, for every $\displaystyle \overline{\xi}\in
R\left(\espro{(-\infty,\alpha]}{\phi(f(a))}\right)$,
$\|\overline{\xi}\|=1$
\[
\alpha\geq\pint{\phi(f(a))\overline{\xi},\;\overline{\xi}}.
\]

Taking this into account, and using Kaplansky's formula we have
that
\begin{align*}
\espro{(\alpha,+\infty)}{f(\phi(a))}&=
\espro{(\alpha,+\infty)}{f(\phi(a))}-\left(\espro{(\alpha,+\infty)}{f(\phi(a))}\wedge
\espro{(-\infty,\alpha]}{\phi(f(a))}\right)\\&\sim
\left(\espro{(\alpha,+\infty)}{f(\phi(a))}\vee
\espro{(-\infty,\alpha]}{\phi(f(a))}\right)-\espro{(-\infty,\alpha]}{\phi(f(a))}\\&\leq
I-\espro{(-\infty,\alpha]}{\phi(f(a))}=\espro{(\alpha,+\infty}{\phi(f(a))}.
\end{align*}
\QED

\medskip
\begin{rems}
$\;$
\begin{itemize}
  \item Note that we can not infer from the above result a similar
  one for monotone concave function because it is not true that $a\leqrt b \Rightarrow -b\leqrt
  -a$. However, using the same arguments a similar result can be proved for monotone concave
  function.
  \item If $0\in I$ and $f(0)\leq 0$ the same result holds for
  positive contractive maps. The proof follows essentially the same lines
  but we have to use Corollary \ref{jensen para comm}, which will be proved
  in the next section, instead of Theorem \ref{jensen para estados}. \EOE
\end{itemize}
\end{rems}

\medskip

In the particular case when $\phi(a)$ and $\phi(f(a))$ are compact
operators and $f(0)=0$ we obtain the following Corollary:

\begin{cor}\label{jensen para compactos}
Let $\cA$ be an unital $C^*$-algebra, $\hil$ a Hilbert space and
$\phi:\cA\rightarrow L(\hil)$ a positive unital map. Then, for
each monotone convex function $f$ defined on $[0,+\infty)$ such
that $f(0)=0$, and for each positive element $a$ of $\cA$ such
that $\phi(a)$ and $\phi(f(a))$ are compact, there exists a
partial isometry $u\in\op$ with initial space
$\overline{R(f(\phi(a)))}$ such that:
\begin{align*}\label{jensen orden truch}
uf(\phi(a))u^*\leq \phi(f(a))  &&\mbox{and}&&
\Big{(}uf(\phi(a))u^* \Big{)}\,
\phi(f(a))=\phi(f(a))\,\Big{(}uf(\phi(a))u^* \Big{)}.
\end{align*}
\end{cor}

The proof of this Corollary is based in the next technical Lemma.

\begin{lem}\label{Lema para meter isometrias}
Let $\hil$ be a Hilbert space. Given $a,b$ positive operators of
$\op$ such that $a$ is compact, $b$ is diagonalizable and $a\leqrt
b$. Then there exists a partial isometry $u\in\op$ with initial space
$\overline{R(a)}$ such that
\begin{align*}
uau^* \leq b&&\mbox{and}&& (uau^*)b=b(uau^*).
\end{align*}
Moreover, if $\hil$ has finite dimension, the same inequality
holds for $a,b\in\opsa$ and the isometry can be changed by a
unitary operator.
\end{lem}
\paragraph{Proof.}
Let $\{\xi_n\}_{n\in \mathbb{N}}$ be an orthonormal basis of
$\overline{R(a)}$ which consists of eigenvectors associated to a
decreasing sequence of positive eigenvalues $\{\la_n\}_{n\in
\mathbb{N}}$ of $a$, counted with multiplicity.

On the other hand, let $\{\eta_m\}_{m\in M}$ be an orthonormal
basis of $\overline{R(b)}$ which consists of eigenvectors
associated to a sequence of positive eigenvalues $\{\mu_m\}_{m\in
M}$ of $b$, counted with multiplicity.

Now, given $\alpha\geq 0$ let $\la_1,\ldots,\la_{n_\alpha}$ be the
eigenvalues of $a$ which are greater than $\alpha$. In a similar
way, let $M(\alpha)=\{m\in M,\, \mu_m > \alpha\}$. If $m_{\alpha}$
denotes be the cardinal of $M(\alpha)$ then notice that $0\leq
m_\alpha\leq \infty$ and if $\alpha\geq \beta\geq 0$ then
$M(\alpha)\subseteq M(\beta)$.

By hypothesis $\espro{(\alpha,+\infty)}{a}$ is equivalent to a
subprojection of $\espro{(\alpha,+\infty)}{b}$. Since
$\espro{(\alpha,+\infty)}{a}$ is the (orthogonal) projection onto
the linear span of $\{\xi_{n}\}_{n=1} ^{n_\alpha}$ and
$\espro{(\alpha,+\infty)}{b}$ is the projection onto the linear
span of $\{\eta_{m}\}_{m\in M(\alpha)}$, then we have that
$n_\alpha\leq m_\alpha$.

\medskip

If we take this into account, we can define an injection
$\psi:\mathbb{N}\rightarrow M$ such that for all $\alpha\geq 0$ it
holds $\psi(k)\in M(\alpha),\, k=1,\ldots, n_{\alpha}$. But then,
defining $u:\overline{R(a)}\rightarrow \overline{R(b)}$ by
\[
u(\xi_n)=\eta_{\psi(n)}
\]
and extending this definition to $\hil$ as zero in $\ker(a)$,  we
obtain that $uau^*\leq b$ and $(uau^*)b=b(uau^*)$, since the set
$(\eta_m)_{m\in M}$ is  a system of eigenvectors for $uau^*$ and $b$, which is complete
for $\overline{R(uau^*)} = R(u)$, and  $\ker ( uau^* )= R(u)^\perp $ is an invariant subspace for  $b$.

Finally, if $\dim\hil=n<\infty$ and $a,b\in\opsa$, let
$\xi_1,\dots, \xi_n$ be a basis of $\hil$ which consists of
eigenvectors of $a$ associated to a decreasing sequence of
eigenvalues $\la_1 , \dots, \la_n$ and let $\eta_1,\dots, \eta_n$
be a basis  of $\hil$ which consists of eigenvectors of $b$
associated to a decreasing sequence of eigenvalues $\mu_1 , \dots
\mu_n$.

Then, if we define $u:\hil\rightarrow\hil$ by mean of
$u(\xi_m)=\eta_m$, $1\le m \le n$, the same argument used before
shows that $uau^*\leq b$ and $(uau^*)b=b(uau^*)$. By construction,
$u$ is unitary.
 \QED

\subsection{Arbitrary convex functions.}\label{Convex functions. Majorization.}

\medskip

As the following example due to J. S. Aujla and F. C. Silva shows,
Theorem \ref{teo de jensen con el orden retruch} may be false if
the function $f$ is not monotone.

\medskip
\begin{exa}
Consider the positive map $\phi:\mathcal{M}_4\rightarrow
\mathcal{M}_2$ given by
\[
\phi\left(\left(%
\begin{array}{cc}
  A_{11} & A_{12} \\
  A_{21} & A_{22} \\
\end{array}%
\right)\right)=\frac{A_{11}+A_{22}}{2}
\]
Take $f(t)=|t|$ and let $A$ be the following matrix
\[
A=\left(%
\begin{array}{rr|cc}
-2&0& 0&0 \\ 0&0& 0&0 \\ \hline 0&0&1&1 \\ 0&0&1&1
\end{array}%
\right)
\]
Then
\begin{align*}
\phi(f(A))=\left(%
\begin{array}{cc}
  3/2 & 1/2 \\
  1/2 & 1/2 \\
\end{array}%
 \right)&&\mbox{and}&& f(\phi(A))=\left(%
\begin{array}{cc}
  1/\sqrt{2} & 0 \\
  0 & 1/\sqrt{2} \\
\end{array}%
\right).
\end{align*}
Easy calculations shows that
$$
\mbox{rank}(\espro{(0.5,+\infty)}{\phi(f(A))})=1<2=
\mbox{rank}(\espro{(0.5,+\infty)}{f(\phi(A))}).
$$
\EOE
\end{exa}

\medskip

Nevertheless, a Jensen's type  inequality holds with respect to
the usual order for every convex function, if the map $\phi$ takes
values in a commutative algebra $\cB$. More generaly, it suffices
to assume that the elements $\phi(f(a))$ and $f(\phi(a))$ of $\cB
$ commute.

\begin{teo}\label{jensen cuando conmutan en la imagen}
Let $\cA$ and $\cB$ be unital \cstars, $\phi:\cA\rightarrow\cB$ a
positive unital map, and $f$ a convex function defined on some
open interval $I$. Given a self-adjoint element $a\in \cA$ whose
spectrum is contained in $I$, if $\phi(a)$ and $\phi(f(a))$
commute then it holds that:
\begin{equation}\label{jensen imagen conmutativa}
  f(\phi(a))\leq \phi(f(a)).
\end{equation}
Moreover, if $0\in I$ and $f(0)\leq 0$ then equation
(\ref{jensen imagen conmutativa}) holds if  $\phi$ is contractive.
\end{teo}

\paragraph{Proof.}

Let $\widehat{\cB}$ denote the abelian C$^*$-subalgebra of $\cB$
generated by $\phi(a)$ and $\phi(f(a))$. On the other hand, let
$\{f_{i}\}_{i\geq 1}$ be the linear functions given by proposition
\ref{propiedad lla}, such that for each $x\in \sigma(a)$
\[
f(x)=\sup_{i\geq 1}f_i(x).
\]

Since $f\geq f_i$ ($i\geq 1$) we have that $f(a)\geq f_i(a)$ and
therefore we obtain that
\begin{equation}\label{f mas grande que cada fi}
\phi(f(a))\geq \phi(f_i(a))=f_i(\phi(a))\, ,
\end{equation}
where the last equality holds because $f_i$ is linear.

As $f_i(\phi(a))$ belongs to the abelian \cstar $\widehat{\cB}$
for every $i\geq 1$, and  $\phi(f(a))$ also belongs to the abelian
\cstar $\widehat{\cB}$ we have that
\[
\phi(f(a))\geq \max_{1\leq i\leq
n}\left[f_i(\phi(a))\right]=\max_{1\leq i\leq
n}\left[f_i\right](\phi(a)).
\]
Now, since $\displaystyle\max_{1\leq i\leq
n}\left[f_i\right]\rightarrow f$ uniformly on compact sets by
Dini's theorem, we obtain
\[
\phi(f(a))\geq f(\phi(a)).
\]

 If $\phi$ is contractive and $f(0)\leq 0$, the
functions $f_i$ also satisfy that $f_i(0)\leq 0$ and we can
replace (\ref{f mas grande que cada fi}) by:
\[
\phi(f(a))\geq \phi(f_i(a))\geq f_i(\phi(a))\, ,
\]
and then we can repeat the same argument to get the desired
inequality. \QED

\medskip

\begin{cor}\label{jensen para comm}
Let $\cA$, $\cB$ be \cstars and suppose that $\cB$ is commutative.
Let $f$ be a convex function defined on some open interval $I$ and
$\phi:\cA\rightarrow \cB$ positive unital. Then
\begin{equation}\label{above equation}
f(\phi(a))\leq \phi(f(a))) \, ,
\end{equation}
for every $a\in\casa$ whose spectrum is contained in $I$.
Moreover, if $0\in I$ and $f(0)\leq 0$ then equation (\ref{above
equation}) holds even if $\phi$ is contractive.
\end{cor}

\bigskip

Let $\cC \subseteq\cB$ be \cstarr s. A \textbf{conditional
expectation} $\cE:\cB\rightarrow \cC$ is a positive $\cC$-linear
projection from $\cB$ onto $\cC$ of norm 1. For instance, states
are conditional expectations. It is well known that a conditional
expectation is a completely positive map. The \textbf{centralizer}
of $\cE$ is the $C^*$-subalgebra of $\cB$ defined by:
\[
\cB^\cE=\{b\in \cB:\ \cE(ba)=\cE(ab),\ \forall a\in \mathcal{B}\}
\]
Notice that, for every $b\in \cB^\cE$, $\cE(b)\in \cZ(\cC)$; where
$\mathcal{Z}(\mathcal{C})=\{c\in \mathcal{C}: \ cb=bc, \ \forall
b\in \mathcal{C}\}$ is the \textbf{center} of $\mathcal{C}$.

\medskip

In \cite{[HanPed]} Hansen and Pedersen considered the vector space
$C_b(X,\cA)$ of all continuous and norm bounded function from a
locally compact Hausdorff space $X$, endowed with a borel measure
$\mu$, to a unital \cstar $\cA$. As it is well known, this space
is a \cstar with the pointwise sum, multiplication, involution and
the norm
\[
\ninf{g}=\sup_X\|g(x)\|.
\]

If the function $x\rightarrow \|g(x)\|$ is integrable, the
function $g$ is called integrable and we can consider the
Bochner's integral
\[
\int_X\; g(x) d\mu(x).
\]
A function $d\in C_b(X,\cA)$ is called \textbf{density} whenever
$d^*d$ is integrable and it holds that
\[
\int_X\; d^*(x)\,d(x)\; d\mu(x)=I.
\]

To each density $d$ there is associated a positive unital map,
$\phi:C_b(X,\cA)\rightarrow \cA$, defined in the following way
\[
\phi(g)=\int_X\; d^*(x)\,g(x)\,d(x)\; d\mu(x).
\]

In this context, given a state $\fii$ of $\cA$ and a convex
function $f$ defined on an open interval $I$, Hansen and Pedersen
proved (\cite{[HanPed]}) the following Jensen's type inequality $$
\varphi\left(f\left(\int_X d^*(X) \, g(x) \, d(x)\,
d\mu(x)\right)\right)\leq \varphi\left(\int_X
d^*(x)\,f(g(x))\,d(x)\,d\mu(x)\right)$$ for each $g\in
C_b(X,\cA)_{sa}$ such that $h=\int_X d^*(x) \, g(x)\,d(x) \,
d\mu(x)\in \cA^\varphi $ and the spectrum of $g(x)$ is contained
in $I$ for every $x\in X$.

Using Stinespring's Theorems, we can rewrite the above result in
the following way:
\begin{equation}\label{Pedersen con stinespring}
\fii(f(\phi(a)))\leq \fii(\phi(f(a)))
\end{equation}
where $\phi:\cA\to\cB$ is a positive unital map between
$C^*$-algebras, and we suppose that $\phi(a)\in\cB^\fii$.

The following Theorem is a slight improvement of
(\ref{Pedersen con stinespring}), which let us connect
Hansen-Pedersen's result with the theory of majorization in finite
factors.

\begin{teo}\label{teo de jensen para trazas}
Let $\cA$ and $\cB$ be unital \cstars and $\phi:\cA\rightarrow\cB$
be a positive unital map. Let $\cE:\cB\rightarrow \cC$ be a
conditional expectation  from $\cB$ onto a $C^*$-subalgebra $\cC$.
Then, for every convex function $f$ defined on some open interval
$I$,

\begin{equation}\label{jensen mayorizacion encubierta}
\cE(g[f(\phi(a))])\leq \cE(g[\phi(f(a))])
\end{equation}
where $a\in\cA_{sa}$ such that $\sigma(a)\subseteq I$ and
$\phi(a)$ belongs to $\cB^\cE$, and $g:J\rightarrow \mathbb{R}$ is
any increasing convex function from some open interval $J$ such
that Im$(f)\subseteq J$ .
\end{teo}

Before starting the proof of this Theorem, we need to prove the
following Lemma.

\begin{lem}\label{lema de la traza}
Let $\cB$ be a $C^*$-algebra, $\fii$ a state defined on $\cB$, and
$b\in \cB^\fii$. Then, there exist a Borel measure $\mu$
defined on $\sigma(b)$ and a positive unital linear map
$\Psi:\cB\rightarrow L^\infty(\sigma(b),\mu)$ such that:
\begin{description}
  \item[i.] $\Psi(f(b))=f$ for every $f\in C(\sigma(b))$.
  \item[ii.]$\displaystyle \fii(x)=\int_{\sigma(b)} \Psi(x)(t)\;d\mu(t)$
  for every $x\in\cB$.
\end{description}
\end{lem}
\paragraph{Proof.}
First of all, note that for every continuous function $g$ defined
on the spectrum of $b$, the map
\[
g\rightarrow \fii(g(b))\, ,
\]
is a bounded linear functional on $C(\sigma(b))$. Therefore, by the
Riesz's representation theorem, there
exists a Borel measure $\mu$ defined on the Borel subsets of
$\sigma(b)$, such that for every continuous function $g$ on
$\sigma(b)$,
\[
\fii(g(b))=\int_{\sigma(b)}g(t)\;d\mu(t).
\]
Now, given $x\in\cB^+$, define the following functional on
$C(\sigma(B))$
\[
g\rightarrow \fii(x g(b)).
\]
Since for every positive element $y$ of $C(b)$,
$
\fii(x y)=\fii(y^{1/2}x y^{1/2})\leq \|x\|\fii(y)\, ,
$
this functional is not only bounded but also
dominated by the functional defined before. So, there exists an
element $h_x$ of $L^{\infty}(\sigma(b),\mu)$ such that, for every
$g \in C(\sigma(b))$,
\[
\fii(x g(b))=\int_{\sigma(b)}g(t)\;h_x(t)\;d\mu(t).
\]
The map $x\mapsto h_x$, extended by linearization, defines a
positive unital linear map $\Psi : \cB \to
L^{\infty}(\sigma(b),\mu)$ which satisfies condition (i), because
\[
\fii(f(b)g(b))=\fii(f g(b))=\int_{\sigma(b)}\;g(t)f(t)\;d\mu(t).
\]
 In order to prove (ii), note that
\[
\fii(x)=\fii(x \ 1(b))=
\int_{\sigma(b)}1\; \Psi(x)(t)\;d\mu(t)=\int_{\sigma(b)}\Psi(x)(t)\;d\mu(t).
\]
\QED

\paragraph{Proof of Theorem \ref{teo de jensen para trazas}.}

Suppose firstly  that $\cE$ is a state. Define
$b=\phi(a)$. Since $b\in \cB^\fii$, by the previous lemma there
exist a Borel measure $\mu$ defined on $\sigma(b)$ and a positive
unital linear map $\Psi:\cB\rightarrow L^\infty(\sigma(b),\mu)$
such that:
\begin{description}
  \item[i.] $\Psi(f(b))=f$ for every $f\in C(\sigma(b))$.
  \item[ii.]$\displaystyle \fii(x)=\int_{\sigma(b)} \Psi(x)(t)\;d\mu(t)$
  for every $x\in\cB$.
\end{description}

\noi
Now, consider the map $\Phi:C(\sigma(a))\rightarrow
L^{\infty}(\sigma(b),\mu)$ defined by
$
\Phi(h)=\Psi(\phi(h(a))).
$
Then, $\Phi$ is bounded, unital and positive. Moreover,
\[
\fii(\phi(h(a)))=\int_{\sigma(b)}\;\Psi(\phi(h(a)))(t)\;d\mu(t)=\int_{\sigma(b)}\Phi(h)(t)\;d\mu(t).
\]
Let $g$ be a increasing convex function. Then, using
Theorem \ref{jensen cuando conmutan en la imagen},
\begin{align*}
\fii(g[f(\phi(a))])&=\fii(g\circ f(b))=\int_{\sigma(b)}g\circ
f(t)\;d\mu(t)=\int_{\sigma(b)}g[f(\Phi(Id))]\;d\mu(t)\\
&\leq\int_{\sigma(b)}g[\Phi(f)(t)]\;d\mu(t)=\int_{\sigma(b)}g[\Psi(\phi(f(a)))(t)]\;d\mu(t)\\
&\leq\int_{\sigma(b)}\Psi(g[\phi(f(a))])(t)\;d\mu(t)=\fii(g[\phi(f(a))]).
\end{align*}
The general case can be reduced to the case already proved by
composing the conditional expectation with the states of the
\cstar $\cC$. \QED

\bigskip
\noi Using Theorem \ref{teo de jensen para trazas}, we can show a
Jensen's type inequality for arbitrary convex functions with
respect to the majorization preorder:

\begin{teo} \label{jensen con mayo}
Let $\cA$ be a unital $C^*$-algebra, $\cB$ a finite factor and
$\phi:\cA\rightarrow\cB$ a positive unital map. Then, for every
convex function $f$, defined on some open interval $I$, and for
every self-adjoint element $a$ of $\cA$ whose spectrum is
contained in $I$,
\begin{equation}\label{jensen con mayorizacion}
f(\phi(a)))\leqm \phi(f(a)).
\end{equation}
\end{teo}
\paragraph{Proof.}
By Theorem \ref{teo de jensen para trazas}, applied to the tracial
state of the finite factor $\cB$,
$$\tr[g(f(\phi(a)))]\leq \tr[g(\phi(f(a)))]$$
for every increasing convex function $g$. By
Proposition \ref{proposicion mayorizacion}, we get
$f(\phi(a)))\leqm \phi(f(a))$. \QED

\begin{rem}
As it happens with Theorems \ref{jensen cuando conmutan en la
imagen} and \ref{teo de jensen con el orden retruch}, in all other
Theorems of this section we can ask $\phi$ to be contractive
instead of unital. In that case the function $f$ has to satisfy
some condition at zero, for example $f(0)\leq 0$. \EOE
\end{rem}

\section{ The multi-variable case.}
In this section we shall be concerned with the restatement, in the
multi-variable context, of several results obtained in section
\ref{Convex functions. Majorization.}. For related results, 
see \cite{[HanPed2]} and \cite{[Ped]}.

\subsection{Multivariated functional calculus}

In what follows we consider the functional calculus for a function
of several variables, so we state a few facts about it in order to
keep the text selfcontained.
 Let
$\mathcal{A}$ be a unital C$^*$-algebra and let $a_{1}, \dots ,a_{n}$
be mutually commuting elements of $\mathcal{A}_{sa}$ i.e., the
self-adjoint part of $\mathcal{A}$. If
$\mathcal{B}=C^*(a_{1}, \dots ,a_{n})$ denotes the unital
C$^*$-subalgebra of $\mathcal{A}$ generated by these elements,
then $\mathcal{B}$ is an abelian $C^*$-algebra. So there exists a
compact Hausdorff space $X$ such that $\mathcal{B}$ is
$^*$-isomorphic to $C(X)$. Actually $X$ is (up to homeomorphism)
the space of characters of $\mathcal{B}$ i.e, the set of
homomorphisms $\gamma:\mathcal{B}\rightarrow \mathbb{C}$, endowed
with the weak*- topology.

\smallskip

Recall that in the case of one operator $a\in \mathcal{A}_{sa}$,
$X$ is homeomorphic to $\sigma(a)$. In general, characters of the
algebra $\mathcal{B}$ are associated in a continuous and injective
way to $n$-tuples $(\lambda_{1}, \dots ,\lambda_{n})\in\prod_{i=1} ^n
\sigma(a_{i})$ by restriction to the abelian C$^*$-subalgebras of
$\mathcal{B}$ generated by the $a_{i}'\, s$. Thus $X$ is
homeomorphic to its image under this map, which we call joint
spectrum and denote $\sigma(a_{1}, \dots ,a_{n})$.

\begin{exa}
Let $\mathcal{A}$ be a unital $C^*$-algebra and $a,\, b\in
\mathcal{A}_{sa}$ such that $ab=ba$. Although the joint spectrum
$\sigma(a,b)$ is a closed subset of the product $\sigma(a)\times
\sigma(b)$ it may be, in general, quite $thin$ , for example
$\sigma(a,b)^0=\emptyset$. Indeed if $b=f(a)$ for a continuous
function $f:\sigma(a)\rightarrow \mathbb{R}$ then is easy to see
that $\sigma(a,b)=\{(x,f(x)),\,x\in \sigma(a)\}=$Graph$(f)$.
Notice that in this case $C^*(a)=C^*(a,b)$ and $\sigma(a)$ is
homeomorphic to $\sigma(a,b)$ in the obvious way. \EOE

\end{exa}

We can now describe the functional calculus in several variables.
Let $f: \sigma(a_{1}, \dots ,a_{n})\rightarrow \mathbb{R}$ be a
continuous function defined on the joint spectrum of the
$a_{i}$'s. Then there exists an element, denoted by
$f(a_{1}, \dots ,a_{n})$, that corresponds to the continuous function
$f$ by the above $^*$-isomorphism. Note that by Tietze's extension
theorem we can consider functions defined on $\prod_{i=1} ^n
\sigma(a_{i})\subseteq \mathbb{R}^n$ without loss of generality.
Therefore the association $f\mapsto f(a_{1}, \dots ,a_{n})$ is a
$^*$-homomorphism from $C(\prod_{i=1} ^n \sigma(a_{i}))$ onto
$\mathcal{B}$, which extends the usual functional calculus of one
variable.

\smallskip

\subsection{Jensen's type inequality in several variables.
Majorization}

\medskip

\begin{fed} \rm
Let $U$ be a convex subset of $\mathbb{R}^n$. A function
$f:U\rightarrow \mathbb{R}$ is called convex if for all $x,\, y\in
U$ and for all $0\leq \lambda \leq 1$ it holds that
$$f(\lambda x + (1-\lambda) y)\leq \lambda f(x)+ (1-\lambda)
f(y).$$
\end{fed}

\smallskip

\begin{pro}\label{lineales}
Let $f$ be a convex function defined on an open convex set
$U\subseteq \mathbb{R}^n$ and let $K\subseteq U$ be compact. Then,
there exists a countable family of linear functions
$\{f_{i}\}_{i\geq 1}$ such that for each $x\in K$ it holds that $$
f(x)=\sup_{i\geq 1} f_{i}(x).$$
\end{pro}

The following results are the multi-variable versions of Theorems
\ref{jensen cuando conmutan en la imagen},  \ref{teo de jensen
para trazas} and  \ref{jensen con mayo} and Lemma \ref {lema de la
traza}. The proofs of those results were chosen in such a way that
they still hold in the multivariable case without substantial
differences.

\begin{teo}\label{gen}
Let $\mathcal{A}$ and $\mathcal{B}$ be unital C$^*$-algebras,
$\phi:\mathcal{A}\rightarrow \mathcal{B}$ be a positive unital map
and $f:U\rightarrow \mathbb{R}$ a convex function. Let
$a_{1}, \dots ,a_{n}\in \mathcal{A}_{sa}$ be mutually commuting, with
$\prod_{i=1} ^n \sigma(a_{i})\subseteq U$ and such that
$\phi(a_{1}), \dots ,\phi(a_{n}), \phi(f(a_{1}, \dots ,a_{n}))$ are also
mutually commuting then
\begin{equation}\label{prop1} f(\phi(a_{1}), \dots ,\phi(a_{n}))\leq
\phi(f(a_{1}, \dots ,a_{n})).
\end{equation}
Moreover, if $\tilde{0} = (0, \dots ,0) \in U$ and $f(\tilde{0})\leq
0$ then equation (\ref{prop1}) holds even if $\phi$ is positive
contractive.
\end{teo}

\medskip

\begin{lem}\label{lema de la traza II}
Let $\cB$ be a $C^*$-algebra, $\fii$ a state defined on $\cB$, and
$b_1, \dots , b_n \in \cB^\fii$ mutually commuting.
Then, there exist a Borel measure $\mu$
defined on $K := \sigma(b_1, \dots , b_n)$ and a positive unital linear map
$\Psi:\cB\rightarrow L^\infty(K,\mu)$ such that:
\begin{description}
  \item[i.] $\Psi(f(b_1, \dots , b_n ))=f$ for every $f\in C(K)$.
  \item[ii.]$\displaystyle \fii(x)=\int_{K} \Psi(x)(t)\;d\mu(t)$
  for every $x\in\cB$.
\end{description}
\end{lem}

\medskip

\begin{teo}

Let $\mathcal{A}$ and $\mathcal{B}$ be unital C$^*$-algebras,
$\phi:\mathcal{A}\rightarrow \mathcal{B}$ be a positive unital
map. Suppose there exists a conditional expectation
$\mathcal{E}:\mathcal{B}\rightarrow \mathcal{C}$, from
$\mathcal{B}$ onto the C$^*$-subalgebra $\mathcal{C}$. Then for
every convex function $f:U\rightarrow \mathbb{R}$ defined on some
open convex set $U\subseteq \mathbb{R}^n$ it holds that
\begin{equation}\label{eq2}
\mathcal{E}(g[f(\phi(a_{1}), \dots ,\phi(a_{n}))])\leq
\mathcal{E}(g[\phi(f(a_{1}, \dots ,a_{n}))]).
\end{equation}
 where
$a_{1}, \dots ,a_{n}\in \mathcal{A}_{sa}$ are mutually commuting
 such that $\prod_{i=1} ^n \sigma(a_{i})\subseteq
U$, $\phi(a_{1}), \dots ,\phi(a_{n})\in \mathcal{B}^\mathcal{E}$
are mutually commuting, and $g:I\rightarrow \mathbb{R}$ is a
convex increasing function defined on some open interval, for
which Im$(f)\subseteq I$.
\end{teo}

\medskip

\begin{teo}
In the conditions of the above theorem let suppose further that
$\mathcal{B}$ is a finite factor, then we have
\begin{equation}\label{jensen con mayorizacion2} f(\phi(a_{1}), \dots ,\phi(a_{n})))\leqm
\phi(f(a_{1}, \dots ,a_{n})).
\end{equation}
\end{teo}

Note that in these Theorems, the requirement that the elements
$a_{1}, \dots ,a_{n}$ (respectively  $\phi(a_{1}), \dots ,\phi(a_{n})$) commute is
needed in order to compute $f(a_{1}, \dots ,a_{n})$ (respectively
$f(\phi(a_{1}), \dots ,\phi(a_{n}))$).  On the other hand the
hypothesis that $\phi(f(a_{1}, \dots ,a_{n}))$ commutes with the
elements $\phi(a_{1}), \dots ,\phi(a_{n})$ corresponds to a technical
reason just as in theorem \ref{jensen cuando conmutan en la
imagen}.

\section{The finite dimensional case.}\label{el caso finito dimensional}

Both the spectral preorder and the majorization have well known
characterizations in the finite dimensional case. The main goal of
this section is to rewrite the already obtained Jensen's
inequalities in terms of this characterizations. Throughout this
section we identify the space $L(\mathbb{C}^n)$ with the space of
complex matrices $\mat$, the real vector space $\opsa$ with the
real vector space $\matsa$ of selfadjoint matrices and the
positive cone $\posop$ with the positive cone of positive matrices
$\matpos$. Given a selfadjoint matrix $A$, by means of
$\lambda_{1}(A), \ldots, \lambda_n(A)$ we denote the eigenvalues
of $A$ counted with multiplicity and arranged in non-increasing
order.

Let us start by recalling the aspect of the spectral preorder and
majorization in $\matsa$.

\begin{num}\label{chamullo}Let $A,B\in\matsa$.
\begin{enumerate}
  \item Using Corollary \ref{Lema para meter isometrias} it is
  easy to see that the following conditions are equivalent
\begin{enumerate}
  \item $A\leqrt B$.
  \item There is an unitary matrix $U$ such that $(UA U^*)B=B(UA
  U^*)$ and $UA U^*\leq B$.
  \item $\avi{i}{A}\leq \avi{i}{B}$ ($1\leq i\leq n$).
\end{enumerate}

  \item Straightforward calculations show that, given a selfadjoint matrix
  $C$,
  the functions $e_C(t)$ considered in the definition of
  majorization have the following form
  \[
e_C(t)=\left\{\begin{array}{l}
  \avi{1}{C}\hspace{0.7cm} \mbox{if}\;\;0\leq t<\frac{1}{n} \\
  \vdots \\
  \avi{n}{C}\hspace{0.7cm} \mbox{if}\;\;\frac{n-1}{n}\leq t<1
\end{array}
\right.
  \]
  Therefore, it holds
\begin{align*}
\int_0^{\alpha}e_{A}(t)\,dt\leq\int_0^{\alpha}
e_{B}(t)\,dt\;\;(\forall\,\alpha\in\mathbb{R}^+)&&\Longleftrightarrow&&
\sum_{i=1} ^k \avi{i}{A} \leq \sum_{i=1} ^k \avi{i}{B}\;\;
(k=1,\ldots , n).
\end{align*}
\end{enumerate}

\end{num}

In the next Proposition, we summarize the different versions of
Jensen's inequality obtained in section 3 using the above
characterizations of the spectral preorder and majorization.

\begin{pro}\label{label largo para que demetrio se enoje}
Let $\cA$ be an unital \cstar and $\phi:\cA\rightarrow \mat$ a
positive unital map. Suppose that $a\in\cA_{sa}$ and that
$f:I\rightarrow \mathbb{R}$ is a function whose domain is an
interval which contains the spectrum of $a$. Then
\begin{description}
  \item[1.] If $f$ is an operator convex function
\[
  \phi(f(a))\leq f(\phi(a)).
\]
  \item[2.] If $f$ is a monotone convex function, for every $i\in\{1,\ldots,n\}$
\[
\avi{i}{f(\phi(a))}\leq \avi{i}{\phi(f(a))}.
\]

  \item[3.] If $f$ is a convex function
\[
\sum_{i=1} ^k \avi{i}{f(\phi(a))} \leq \sum_{i=1} ^k
\avi{i}{\phi(f(a))}.
\]

Moreover, if $0\in I$ and $f(0)\leq 0$ the above inequalities also
hold for contractive positive maps.
\end{description}
\end{pro}

\begin{exa}\label{multiplicador Schur}
Given two $n\times n$ matrices $A=(a_{ij})$ and $B=(b_{ij})$, we
denote $A\circ B=(a_{ij}b_{ij})$ their Schur's product. It is a
well known fact that the map $A\mapsto A\circ B$ is completely
positive for each positive matrix $B$, and if we further assume
that $I\circ B=I$ then it is also unital. So the above
inequalities can be rewritten taking $\phi:\mat\rightarrow\mat$
given by $\phi(A)=A\circ B$ where $B$ satisfies the mentioned
properties. \EOE
\end{exa}

\begin{exa}
Another example can be obtained by taking
$\phi:\mathcal{M}_m\rightarrow \mat$ given by
\[
\phi(A)=\sum_{i=1} ^r W_i ^* A\,W_i
\]
where $W_1,\ldots, W_r$ are $m\times n$ matrices such that
$\displaystyle\sum_{i=1} ^r W_i ^* \,W_i=I$. These maps are
completely positive by Choi's Theorem and they are unital as a
result of the condition over the rectangular matrices $W_i$. \EOE
\end{exa}

\begin{exa}
As a final application let us consider for each $\alpha\in(0,1)$
the positive unital map $\phi_\alpha: \mathcal{M}_{2n} \rightarrow \mat$
defined by
\[
\phi\left(\begin{pmatrix}
  A & B \\
  C & D
\end{pmatrix}\right)=\alpha A+(1-\alpha)D.
\]
Using these maps in Proposition \ref{label largo para que demetrio
se enoje} and taking diagonal block matrices we get that for every
monotone convex functions $f$
\begin{align*}
&\avi{i}{f(\alpha A+(1-\alpha)D)}\leq \avi{i}{\alpha
f(A)+(1-\alpha)f(D)} \hspace{2cm}(i=1,\ldots, n)\\ \intertext{ and
for general convex functions $f$} &\sum_{i=1}^r \avi{i}{f(\alpha
A+(1-\alpha)D)}\leq \sum_{i=1}^r \avi{i}{\alpha
f(A)+(1-\alpha)f(D)} \hspace{0.7cm}(r=1,\ldots, n)
\end{align*}
where we are assuming that $A$ and $D$ are selfadjoint matrices
whose spectra are contained in the domain of $f$. This
inequalities were proved by J. S. Aujla and F. C. Silva in
\cite{[AujSil]}. \EOE
\end{exa}

\section{Some inequalities.}\label{some inequalities}

In this section we deduce some inequalities from the results
developed in previous sections by choosing particular (operator)
convex functions and positive maps. See \cite{[BatCh]} and \cite{[Lieb]} for related inequalities.

First of all, let us recall that given a locally compact Hausdorff
space $X$ and a unital \cstar $\cA$, $C_b(X,\cA)$ is the space of
functions $g:X\rightarrow \cA$ which are norm continuous and norm
bounded.

When a Borel measure $\mu$ defined on the Borel subsets of $X$,
for those functions\newline $g\in C_b(X,\cA)$ such that
$x\rightarrow \|g(x)\|$ is integrable, we can consider the Bochner
integral
\[
\int_X\; g(x) d\mu(x)
\]

In particular, we shall say that $d\in C_b(X,\cA)$ is a
\textbf{density} whenever it holds that
\[
\int_X\; d^*(x)\,d(x)\; d\mu(x)=I.
\]

Also recall that to each density $d$ there is associated a
positive unital map, $\phi:C_b(X,\cA)\rightarrow \cA$, defined in
the following way
\[
\phi(g)=\int_X\; d^*(x)\,g(x)\,d(x)\; d\mu(x).
\]

We can now state the following non-commutative version of the
Information inequality.

\medskip

\begin{pro}[Information inequality]\label{information inequality}
Let $X$ be a locally compact Hausdorff space, $\mu$ a Borel
measure on $X$, $\cA$ a unital $C^*$-algebra, and $\;a,b\in
C_b(X,\cA)$ densities such that $a(x),b(x)\in \inversible(\cA)$
for every $x\in X$. Then, it holds that
\[
\int_X\; a^*\,\log(a^{*-1}\,b^*\,b\, a^{-1})\,a\; d\mu\leq 0.
\]
\end{pro}
\paragraph{Proof.}
By Theorem \ref{jensen pesadazo}, and using the fact that the
function $\log$ is operator concave we obtain
\begin{align*}
\int_X\; a^*\,\log(a^{*-1}\,b^*\,b\, a^{-1})\,a\; d\mu &\leq
\log\left[\int_X\; a^*\,a^{*-1}\,b^*\,b\, a^{-1}\,a\;
d\mu\right]\\ &= \log\left[\int_X\; b^*\,b\; d\mu\right]\leq
\log(I)=0.
\end{align*}
\QED

\medskip
Returning to the general case, the following is a non-commutative version of
Liapounov's inequality.  

\begin{pro}[Liapunov's inequality]\label{liapunov}
Let $\cA$ a unital \cstar and $\hil$ a Hilbert space. Given a
positive contractive map $\phi:\cA\rightarrow \op$ and $1\leq
r\leq s$ it holds that
\[
\phi(a^{r})^{1/r}\leq \phi(a^{s})^{1/s}\hspace{1cm} \forall\;
a\in\cA^+.
\]
\end{pro}
\paragraph{Proof.}
Let $b=a^r$ and $t=s/r$. Then, as the function $f(x)=x^{1/t}$ is
operator concave
\[
\phi(b^t)^{1/t}\geq\phi(b).
\]
So we obtain
\[
\phi(a^s)^{r/s}\geq\phi(a^r)\, ,
\]
and using that the function $f(x)=x^{1/r}$ is operator monotone it
holds that
\[
\phi(a^s)^{1/s}\geq\phi(a^r)^{1/r}.
\]
\QED

Finally let us prove the following H\={o}lder's type inequality
whose matrix version were proved by T. Ando and F. Hiai in
$\cite{[AndoHiai]}$
\begin{pro}
Let $\cA$ be a finite factor and $1<p,q<\infty$ such that
$1/p+1/q=1$. If $c,d\in \cA^+$ and $c^q+d^q\leq I$, then
\[
\tr(ca+db)\leq \Big{(}\tr(a^p+b^p)\Big{)}^{1/p}
\]
for every $a,b\in \cA^+$
\end{pro}
\paragraph{Proof.}
When $c,d\geq 0$ and $c^q+d^q\leq I$, if $c_1=(I-d^q)^{1/q}$, then
$c_1^q+d^q=I$ and $\tr(ca)\leq \tr(c_1a)$ because $c\leq c_1$.
Hence we may assume that $c^q+d^q=I$ so that $c$, $d$ commute.

On the other hand, by Lemma \ref{lema de la traza}, there exist a
measure of probability $\mu$ defined in $\spec{c}$ and a positive
unital map $\Psi:\cA\rightarrow L^\infty(\spec{c},\mu)$ such that
\begin{description}
  \item[i.] $\Psi(f(c))=f$ for every $f\in C(\sigma(c))$.
  \item[ii.]$\displaystyle \tr(x)=\int_{\sigma(c)} \Psi(x)(t)\;d\mu(t)$
  for every $x\in\cB$.
\end{description}
Then
\begin{align*}
\tr(ca+db)&=\int_{\sigma(c)}
t\Psi(a)(t)+(1-t^q)^{1/q}\Psi(b)(t)\;d\mu(t)\\\intertext{using
H\={o}lder inequality for each $t\in \spec{c}$} &=
\int_{\sigma(c)}
\Big{(}t^q+(1-t^q)\Big{)}^{1/q}\Big{(}\Psi(a)^p(t)+\Psi(b)^p(t)\Big{)}^{1/p}\;d\mu(t)\\&=
\int_{\sigma(c)}
\Big{(}\Psi(a)^p(t)+\Psi(b)^p(t)\Big{)}^{1/p}\;d\mu(t)\\
&\leq \left(\int_{\sigma(c)}
\Psi(a^p)(t)+\Psi(b^p)(t)\;d\mu(t)\right)^{1/p}\\&=\Big{(}\tr(a^p+b^p)\Big{)}^{1/p}.
\end{align*}
where the inequality follows from the standard version of Jensen's
inequality for integrals and Theorem \ref{jensen cuando conmutan
en la imagen}. \QED

\subsubsection*{Aknowledgements:} Partially supported by
Universidad de La Plata (UNLP 11 X350) and ANPCYT (PICT 03-09521)


\vglue1truecm

\noindent{\bf Jorge A. Antezana, Pedro G. Massey and Demetrio Stojanoff},

\noi Depto. de Matem\'atica, FCE-UNLP,\\  Calles 50 y 115 (1900), La Plata,  Argentina

\noindent{e-mail: \\antezana@mate.unlp.edu.ar\\ massey@mate.unlp.edu.ar \\ demetrio@mate.unlp.edu.ar}


\begin{thebibliography}{XXXXXX}

\bibitem{[AndoHiai]} Tsuyoshi Ando and Fumio Hiai, H\={o}lder type
inequalities for matrices, Math. Inequalities and Appl. 1 (1998)
1-30.
\bibitem{[AujSil]} Jaspal Singh Aujla and Fernando C. Silva, Weak
majorization inequalities and convex functions, Linear Algebra
Appl. 369 (2003), 217-233.

\bibitem{[BenSher]} Julius Bendat and Seymour Sherman, Monotone
and convex operator functions, Transactions of the American
Mathematical Society 79 (1955), 58-71
\bibitem{[Bat]} Rajendra Bhatia, Matrix Analysis,
Berlin-Heildelberg-New York, Springer 1997.
\bibitem{[BatCh]} Rajendra Bhatia and Chandler Davis, More Operator
Versions of the Schwarz inequality, Communication in Mathematical
Physics, Springer Verlag 215 (2000), 239-244.

\bibitem{[BrKo]} Lawrence G. Brown and Hideki Kosaki, Jensen's
inequality in semi-finite von Neumann algebras, Journal of
Operator Theory 23 (1990), 3-19.


\bibitem{[HanPed]} Frank Hansen and Gert K. Pedersen, Jensen's
operator inequality, Bulletin of London Mathematical Society 35
(2003), 553-564.
\bibitem{[HanPed2]} Frank Hansen and Gert K. Pedersen, 
Jensen's trace inequality in several variables,
Internat. J. Math. 14 (2003), 667-681.

\bibitem{[Ka]} Eizaburo Kamei, Majorization in finite factors,
Math. Japonica 28, No. 4 (1983), 495-499.

\bibitem{[Lieb]} E. H. Lieb and M. B. Ruskai, Some operator
inequalities of the Schwarz type, Advances in Math. 12 (1974),
269-273.



\bibitem{[Pau]} Vern I. Paulsen, Completely bounded maps and dilations,
Notes in Mathematics Series 146, Pitman, New York, 1986.

\bibitem{[Ped]} G. K. Pedersen, Convex trace functions of several variables on 
$C\sp *$-algebras. J. Operator Theory 50 (2003),  157-167.


\bibitem{[Hosh]} Allan M. Sinclair and Roger R. Smith, Hoshschild
Cohomology of von Newmann Algebras, London Math. Soc. Lectures
Notes Series 203, Cambridge Univ. Press, 1995.





\end{thebibliography}
\end{document}